\documentclass[12pt]{article}
\usepackage{amssymb}
\usepackage{latexsym,bm}
\usepackage{graphicx}
\usepackage{amsmath}
\usepackage{mathrsfs}
\usepackage{mathrsfs,amscd,amssymb,amsthm,amsmath,bm,graphicx,psfrag,subfigure,url,xcolor}
\usepackage{tikz}

\usepackage{epstopdf}

\date{}
\newcounter{mathitem}
{\begin{list}{{$(\roman{mathitem})$}}{
\setcounter{mathitem}{0}
\usecounter{mathitem}
\setlength{\topsep}{0pt plus 2pt minus 0pt}
\setlength{\parskip}{0pt plus 2pt minus 0pt}
\setlength{\partopsep}{0pt plus 2pt minus 0pt}
\setlength{\parsep}{0pt plus 2pt minus 0pt}
\setlength{\leftmargin}{35pt}
\setlength{\itemsep}{0pt plus 2pt minus 0pt}}}
{\end{list}}
	
\begin{document}
\title{Hamiltonian graphs with prescribed minimum degree and no near-spanning cycles}
\author{Xingzhi Zhan\thanks{Department of Mathematics, Key Laboratory of MEA (Ministry of Education)
and Shanghai Key Laboratory of PMMP, East China Normal University, Shanghai 200241, China}
\thanks{E-mail address: zhan@math.ecnu.edu.cn}}
\maketitle

\begin{abstract}
In 1984, Roland H\"{a}ggkvist posed the problem of constructing Hamiltonian graphs of order $n$ with large minimum degree and no $(n-2)$-cycle.
He remarked that he did not know of such a graph with minimum degree at least three. We solve this problem by proving the following two results.
(1) For every integer $d\ge 3$ and every integer $n\ge 15d-14,$ there exists a Hamiltonian graph of order $n$ and minimum degree $d$ that contains no
$(n-2)$-cycle. (2) For every integer $d\ge 3,$ every positive integer $k,$ and every integer $n\ge (k+1)[(d-1)(k+3)+1],$ there exists a Hamiltonian graph of order $n$ and minimum degree $d$ that contains no $(n-s)$-cycle for any $s\in\{1,2,\dots,k\}.$ The proofs are constructive. We also pose several open problems.
\end{abstract}


\vskip 8mm

\section{Introduction}

We consider finite simple graphs and use standard terminology and notation from [1] and [8]. The {\it order} of a graph is its number of vertices, and the
{\it size} is its number of edges. We denote by $V(G),$ $E(G),$ and $\delta (G)$  the vertex set, edge set, and minimum degree of a graph $G,$ respectively.
The order of $G$ is denoted by $|G|.$ A {\it $k$-cycle} is a cycle of length $k.$

The study of cycle lengths in Hamiltonian graphs is an interesting area of research that is not yet fully understood [2--7].
For every integer $n\ge 10,$ it is easy to construct a Hamiltonian graph of order $n$ with minimum degree at least three that contains no $(n-1)$-cycle.
The analogous problem for $(n-2)$-cycles is much less straightforward. In 1984, Roland H\"{a}ggkvist [5] posed the following problem.

{\bf Problem 1.} Construct Hamiltonian graphs of order $n$ with large minimum degree and no $(n-2)$-cycle.

He remarked that he did not know of such a graph with minimum degree at least three. We solve this problem by proving the following two results.
The proofs are constructive.

{\bf Theorem 1.} {\it For every integer $d\ge 3$ and every integer $n\ge 15d-14,$ there exists a Hamiltonian graph of order $n$ and minimum degree $d$ that contains no
$(n-2)$-cycle.}

{\bf Theorem 2.} {\it For every integer $d\ge 3,$ every positive integer $k,$ and every integer $n\ge (k+1)[(d-1)(k+3)+1],$ there exists a Hamiltonian graph of order $n$ and minimum degree $d$ that contains no $(n-s)$-cycle for any $s\in\{1,2,\dots,k\}.$}

In posing Problem 1, H\"{a}ggkvist added, \textquotedblleft It is suspected that such graphs exist with minimum degree $cn$ for some $c>0$ (the constant $c$ is necessarily small).\textquotedblright Theorem 1 confirms this suspicion; one may take $c=1/15.$

The degree of a vertex $v$ in a graph $G$ is denoted by ${\rm deg}_G(v),$ and for graphs $G_1, G_2, \dots, G_p$ we write $G_1+G_2+\dots+ G_p$ for the pairwise vertex-disjoint union of these $p$ graphs.

In Section 2, we prove Theorems 1 and 2, and in Section 3, we pose several related open problems.

\section{Proofs of the main results}

{\bf Proof of Theorem 1.}
We will construct a Hamiltonian graph $G=G(d,n)$ of order $n$ with $\delta(G)=d$ that contains no $(n-2)$-cycle.
Put
\[
 h=d-1
 \qquad\text{and}\qquad
 m=n-(15d-14).
\]
Thus $h\ge2$ and $m\ge0$.

\medskip
\noindent\textit{Construction of the vertex set.}
Let
\[
 S=\{x_1,x_2,\ldots,x_{3h}\}
\]
be an independent set of $G$, partitioned into three classes:
\[
 \begin{aligned}
  S_1&=\{x_1,x_2,\ldots,x_h\},\\
  S_2&=\{x_{h+1},x_{h+2},\ldots,x_{2h}\},\\
  S_3&=\{x_{2h+1},x_{2h+2},\ldots,x_{3h}\}.
 \end{aligned}
\]
Class subscripts are cyclic: $S_4$ means $S_1$ and $S_5$ means
$S_2$.  We also put $x_{3h+1}=x_1$.

For $1\le i\le3h$, define
\[
 c_i=
 \begin{cases}
  1,&\text{if }1\le i\le h,\\
  2,&\text{if }h+1\le i\le2h,\\
  3,&\text{if }2h+1\le i\le3h.
 \end{cases}
\]
Thus $x_i\in S_{c_i}$.  Put $c_{3h+1}=1$.  We have
$c_i=c_{i+1}$ except at the three transition indices $h,\,\,2h,\,\,3h.$
Define the positive integers $q_i$ by
\[
 q_i=
 \begin{cases}
  4+m,&\text{if }\,\, i=1,\\
  4,&\text{if }\,\, i\notin\{1,h,2h,3h\},\\
  5,&\text{if }\,\, i\in\{h,2h\},\\
  3,&\text{if }\,\, i=3h.
 \end{cases}
\]
For each $i\in\{1,2,\ldots,3h\}$, introduce a path
\[
 P_i=p_{i,1}p_{i,2}\cdots p_{i,q_i}.
\]
These paths are pairwise vertex-disjoint and are disjoint from $S$. Let $V(G)=S\cup_{i=1}^{3h} V(P_i).$

\medskip
\noindent\textit{Construction of the edge set.}
Include all path edges of every $P_i$.  We now define all edges
between the paths and $S$.

Suppose first that $c_i=c_{i+1}=j$.  Join each of
$p_{i,1}$ and $p_{i,q_i}$ to every vertex of $S_j$, and join every
internal vertex of $P_i$ to every vertex of $S_{j+1}$.

Next, suppose that $i\in\{h,2h\}$ and $c_i=j$.  Then
$c_{i+1}=j+1$ and $q_i=5$.  Join $p_{i,1}$ to every vertex of
$S_j$; join $p_{i,2}$ and $p_{i,5}$ to every vertex of $S_{j+1}$;
and join $p_{i,3}$ and $p_{i,4}$ to every vertex of $S_{j+2}$.

Finally, $P_{3h}$ is the special transition path from $S_3$ to
$S_1$ and has order $3$.  Join $p_{3h,1}$ to every vertex of $S_3$,
and join $p_{3h,2}$ and $p_{3h,3}$ to every vertex of $S_1$.

There are no other edges.  In particular,
\[
 G-S=P_1+P_2+\cdots+P_{3h}.
\]

For illustration, consider the instance $d=3$ and $n=32$.  Here
$h=2$ and $m=1$, and hence
\[
 (q_1,q_2,q_3,q_4,q_5,q_6)=(5,5,4,5,4,3).
\]
The resulting graph $G(3,32)$ is shown in Figure 1.

\begin{figure}[ht]
\centering
\includegraphics[width=0.7\textwidth]{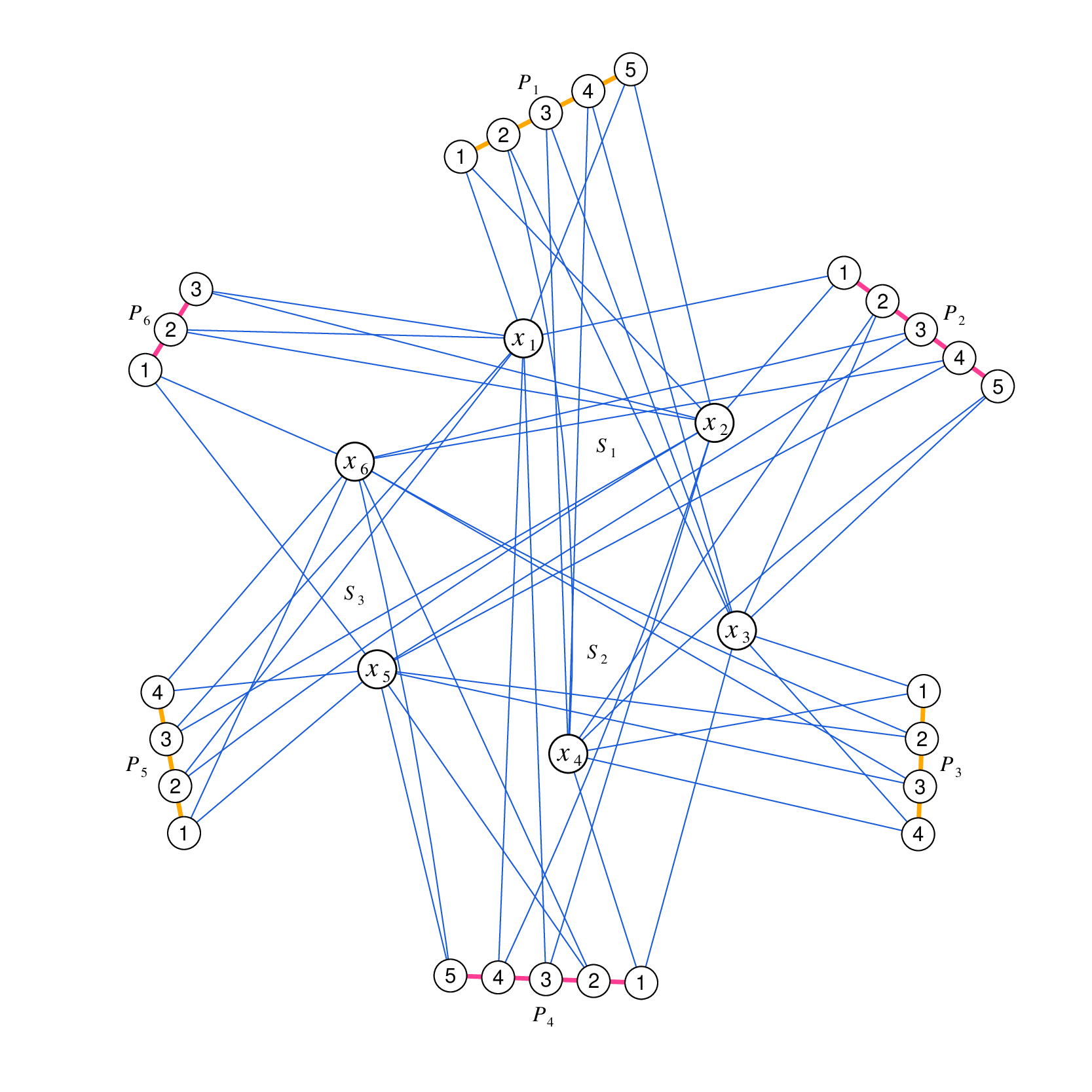}
\caption{The graph $G(3,32)$}
\end{figure}

\medskip
\noindent\textit{Verification of the order.}
There are $3h$ vertices in $S$.  There are $3(h-1)$ indices $i$
satisfying $c_i=c_{i+1}$, and the corresponding paths $P_i$ have
total order
\[
 3(h-1)\cdot4+m=12h-12+m.
\]
The two ordinary transition paths have total order $10$, while the
special transition path has order $3$.  Consequently,
\begin{align*}
 |G|
 &=3h+(12h-12+m)+10+3\\
 &=15h+1+m\\
 &=15(d-1)+1+m\\
 &=15d-14+m=n.
\end{align*}

\medskip
\noindent\textit{Verification that $\delta(G)=d$.}
Every path vertex is adjacent to all $h$ vertices of exactly one
class $S_j$.  A path endpoint has one neighbor on its path, so its
degree is $1+h=d.$ An internal path vertex has two neighbors on its path, so its degree
is $2+h=d+1.$

It remains to check the $S$-vertices.  Fix $v\in S_j$.  For each
$x_i\in S_j$, the left endpoint of $P_i$ is adjacent to $v$.
Moreover, the right endpoint of the path preceding $P_i$ in the
cyclic list is adjacent to $v$.  These are two distinct path endpoints
for each of the $h$ choices of $x_i$.  Therefore
\[
 {\rm deg}_G(v)\ge2h=2(d-1)\ge d,
\]
where the last inequality follows from $d\ge3$.  Thus every vertex
has degree at least $d$, while every path endpoint has degree exactly
$d$.  Hence $\delta(G)=d.$

\medskip
\noindent\textit{Verification of Hamiltonicity.}
For every $i$, the vertex $x_i$ is adjacent to $p_{i,1}$, and
$x_{i+1}$ is adjacent to $p_{i,q_i}$.  Therefore
\[
 x_1P_1x_2P_2\cdots x_{3h}P_{3h}x_1
\]
is a cycle.  It contains every vertex of $S$ and every vertex of
every path exactly once.  Thus it is a Hamilton cycle of $G$.

\medskip
\noindent\textit{Exclusion of an $(n-2)$-cycle.}
Suppose, for a contradiction, that $C$ is a cycle of length $n-2$.
Exactly two vertices of $G$ are omitted from $C$.  Every $P_i$ has at
least three vertices, so $C$ cannot omit an entire path.  Hence $C$
meets all $3h$ paths.

Let $u=|V(C)\cap S|$.  We have $u>0$, because $G-S$ is a disjoint
union of paths and therefore contains no cycle.  Since $S$ is
independent, no two vertices of $V(C)\cap S$ are consecutive on $C$.
Deleting these $u$ vertices from $C$ therefore leaves exactly $u$
non-null path segments.  Each segment is contained in one component
$P_i$ of $G-S$.  Since $C$ meets all $3h$ paths, at least $3h$
segments are required.  Thus $u\ge3h$.  Since $u\le|S|=3h$, we have
$u=3h$.

Since $u=3h=|S|$, the cycle $C$ contains every vertex of $S$.
Deleting $S$ from $C$ leaves exactly $3h$ non-null path segments,
each contained in one of the $3h$ components $P_1,\ldots,P_{3h}$ of $G-S$.  Moreover, $C$ meets every $P_i$, so
each $P_i$ contains at least one such segment.  It follows that each
$P_i$ contains exactly one segment.  This segment is a subpath of
$P_i$, possibly consisting of a single vertex, and therefore consists
of consecutive vertices of $P_i$.  Write the subpath retained from $P_i$ as
\[
 R_i=p_{i,a_i+1}p_{i,a_i+2}\cdots p_{i,q_i-b_i},
 \qquad a_i,b_i\ge0.
\]
Here the omitted \emph{prefix} is the initial segment
$p_{i,1},\ldots,p_{i,a_i}$, and the omitted \emph{suffix} is the
terminal segment
$p_{i,q_i-b_i+1},\ldots,p_{i,q_i}$.  The omitted prefix or suffix is empty when its
corresponding parameter is $0$.  Since all $S$-vertices belong to
$C$, the two omitted vertices occur in these prefixes and suffixes.
Hence
\[
 \sum_{i=1}^{3h}(a_i+b_i)=2. \eqno (1)
\]
Since $R_i$ is a component of $C-S$ and $C$ is a cycle, exactly two edges of $C$ have
one endpoint in $V(R_i)$ and the other endpoint in $S$.  If $R_i$
has at least two vertices, one of these edges is incident with the
left endpoint of $R_i$, and the other is incident with its right
endpoint.  We call these edges the \emph{left token} and the
\emph{right token} of $R_i$, respectively.  If $R_i$ consists of a
single vertex, the two edges are still distinct; we designate one as
the left token and the other as the right token. Thus a token is an
edge of $C$ joining the retained subpath $R_i$ to a vertex of $S$.
The \emph{actual class} of a token is the class $S_j$ containing its
endpoint in $S$. We assign baseline classes to the two tokens of $R_i$ according to
the two endpoints of the original path $P_i$.  The baseline class of the
left token of $R_i$ is $S_{c_i}$, and the baseline class of its right
token is $S_{c_{i+1}}$.

Each class is the baseline class of $h$ left tokens and $h$ right
tokens; hence it is the baseline class of exactly $2h$ tokens.
There are also exactly $2h$ actual tokens in each class, because $C$
contains all $h$ vertices of that class and has degree $2$ at every
one of them.

Orient the cyclic ordering of the three classes as
\[
 S_1\longrightarrow S_2\longrightarrow S_3\longrightarrow S_1.
\]
For each token, define its \emph{forward movement} to be the number
of directed steps, equal to $0$, $1$, or $2$, required to go from
its baseline class to its actual class in this cyclic ordering.
Thus the movement is $0$ when the two classes coincide.
By (1), every $a_i$ and $b_i$ is at most $2$.

Consider a path $P_i$ satisfying $c_i=c_{i+1}=j$, and let
\[
 R_i=p_{i,a_i+1}p_{i,a_i+2}\cdots p_{i,q_i-b_i}
\]
be its retained subpath.  If $a_i=0$, then the left endpoint of
$R_i$ is $p_{i,1}$, which is adjacent to vertices of $S_j$.
Therefore, the left token has movement $0$.  If $a_i>0$, then the
left endpoint $p_{i,a_i+1}$ of $R_i$ is an internal vertex of $P_i$
and is adjacent to vertices of $S_{j+1}$.  Hence the left token has
movement $1$, which is at most $a_i$.  Similarly, the right token
has movement $0$ when $b_i=0$ and movement $1$ when $b_i>0$.
Consequently, the movements of the left and right tokens are at
most $a_i$ and $b_i$, respectively.

For an ordinary transition path, the attachment classes of its five
successive vertices are
\[
 S_j,\ S_{j+1},\ S_{j+2},\ S_{j+2},\ S_{j+1}.
\]
At the left end, the movements corresponding to $a_i=0,1,2$ are
$0,1,2$.  At the right end, the movement is $0$ when $b_i=0$ and is
$1$ when $b_i=1$ or $2$.  Again, the movement at an end is at most
the number of omitted vertices there.

For the special transition path, the attachment classes of its three
successive vertices are
\[
 S_3,\ S_1,\ S_1.
\]
At its left end, the movement is $0$ when $a_{3h}=0$ and is $1$ when
$a_{3h}=1$ or $2$.  At its right end, the movement is $0$ when
$b_{3h}=0$ or $1$, and is $2$ when $b_{3h}=2$.  Once again, movement
is at most the number of omitted vertices at that end.

Let $D$ be the sum of all token movements.  The preceding analysis
and (1) give
\[
 0\le D\le2. \eqno (2)
\]
In fact, $D>0$.  The only nonzero end-trimming that produces
zero movement is the omission of exactly one vertex from the right
end of the special path.  If this happens, the second omitted vertex
produces positive movement.  If two vertices are omitted from that
right end, the right token itself has movement $2$.  Therefore two
omissions cannot produce total movement zero.

For each token, record the directed steps traversed when going from
its baseline class to its actual class along $S_1\longrightarrow S_2\longrightarrow S_3\longrightarrow S_1.$
A traversal of one directed step $ S_j\longrightarrow S_{j+1} $ is called a \emph{unit movement}, where class subscripts are read
cyclically, so that $S_4=S_1$.  Thus a token with movement $0$
produces no unit movement; a token with movement $1$ produces one
unit movement; and a token with movement $2$ produces two consecutive
unit movements.  For $j\in\{1,2,3\}$, let $M_j$ denote the total
number, counted with multiplicity, of unit movements across the
directed step $S_j\longrightarrow S_{j+1}.$

Replace each movement of length $2$ by two consecutive unit
movements. For each class, the number of tokens having that actual
class equals the number having that baseline class. Hence the number
of unit movements entering the class equals the number leaving it,
and therefore
\[
 M_1=M_2=M_3.
\]
Since $D>0$, their common value is positive.  Consequently,
\[
 D=M_1+M_2+M_3\ge3,
\]
contradicting (2).  Thus $G$ contains no cycle of length $n-2$. This completes the proof. \hfill$\Box$

The bound in Theorem 1 gives $n\ge31$ when $d=3$. There is a smaller graph of order $19$ and size $33$ having the three required
properties. See Figure 2.
\begin{figure}[ht]
\centering
\includegraphics[width=0.5\textwidth]{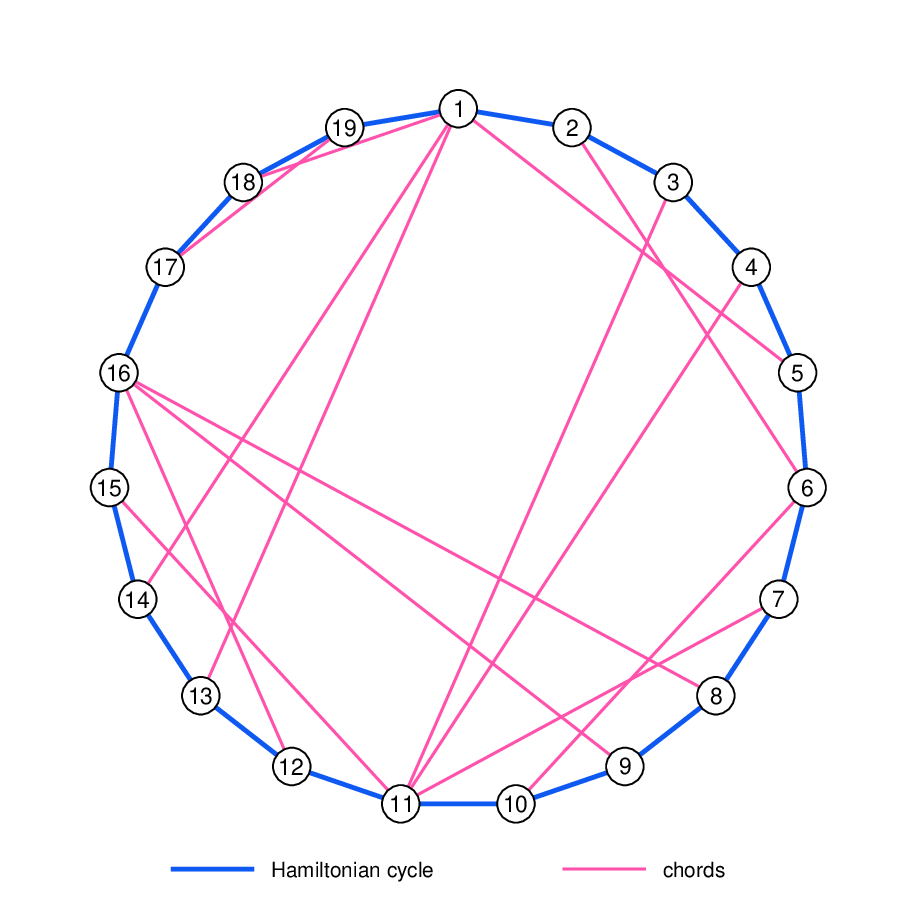}
\caption{A Hamiltonian graph of order $19$ and minimum degree $3$ with no $17$-cycle}
\end{figure}
An exact, exhaustive computation shows that the graph in Figure 2 has the smallest possible order and size among all Hamiltonian graphs of order $n\ge 5$ with minimum degree at least $3$
and no $(n-2)$-cycle. But it is not the unique graph attaining the minimum order and minimum size.

Next, we turn to the proof of Theorem 2.

{\bf Proof of Theorem 2.}
We will construct a graph $G=G(d,k,n)$ of order $n$ such that $\delta(G)=d$ and $G$ contains no cycle of length $n-s$ for any
$s\in\{1,2,\ldots,k\}$.
Put
\[
 f=k+1\qquad\text{and}\qquad h=d-1.
\]
Thus $f\ge2$ and $h\ge2$.

\medskip
\noindent\textit{Construction of the vertex set.}
Let
\[
 T=\{x_1,x_2,\ldots,x_{fh}\}
\]
be an independent set of $G$. For $1\le j\le f$, define
\[
 T_j=\{x_{(j-1)h+1},x_{(j-1)h+2},\ldots,x_{jh}\}.
\]
Thus $T_1,T_2,\ldots,T_f$ partition $T$; each class has $h$
vertices, and $|T|=fh$. Class subscripts and class labels are
interpreted cyclically modulo $f$, with representatives in
$\{1,2,\ldots,f\}$. In particular, $T_{f+1}=T_1$, $T_{f+2}=T_2$,
and the cyclic successor of the label $f$ is $1$.
Regard
\[
 x_1,x_2,\ldots,x_{fh},x_1
\]
as a cyclic listing, and set $x_{fh+1}=x_1$.  Define the class labels
$c_i$ explicitly by
\[
 c_i=j
 \quad\text{if}\quad
 (j-1)h+1\le i\le jh
 \qquad(1\le j\le f),
\]
and set $c_{fh+1}=c_1=1$.  Thus $x_i\in T_{c_i}$ for
$1\le i\le fh+1$.

If $x_i$ and $x_{i+1}$ occur consecutively within the same class,
then $c_i=c_{i+1}$.  There are $h-1$ such consecutive pairs in each
of the $f$ classes, and hence there are exactly $f(h-1)$ indices $i$
for which $c_i=c_{i+1}.$ At the end of each class, the cyclic sequence passes to the next
class. These are the $f$ indices $h,2h,\ldots,fh$. At each such
index, $c_{i+1}$ is the cyclic successor of $c_i$.
We call $h,2h,\ldots,fh$ the \emph{transition indices}.  Thus an
index $i$ is a transition index precisely when the cyclic sequence
passes from $T_{c_i}$ to the next class $T_{c_{i+1}}$; every other
index satisfies $c_i=c_{i+1}$.

Let
\[
 m=n-f\bigl(h(k+3)+1\bigr).
\]
The hypothesis on $n$ and the equalities $f=k+1$ and $h=d-1$ imply
that $m\ge0$.  Since $h\ge2$, the vertices $x_1$ and $x_2$ belong to
the same class $T_1$, so $c_1=c_2=1$.  For each
$i\in\{1,2,\ldots,fh\}$, define
\[
 q_i=
 \begin{cases}
  k+2+m,&\text{if }i=1,\\
  k+2,&\text{if }i\ne1\text{ and }c_i=c_{i+1},\\
  k+3,&\text{if }i\text{ is a transition index}.
 \end{cases}
\]
Having defined $q_i$, introduce the path
\[
 P_i=p_{i,1}p_{i,2}\cdots p_{i,q_i}
 \qquad(1\le i\le fh).
\]
Thus every path has at least $k+2$ vertices, and every transition
path has exactly $k+3$ vertices. These paths are pairwise
vertex-disjoint and are disjoint from $T$.
Let $V(G)=T\cup_{i=1}^{fh} V(P_i).$

\medskip
\noindent\textit{Construction of the edge set.} First include all path edges of every $P_i$.
We next specify all edges between the paths and $T$.  First suppose
that $c_i=c_{i+1}=j$.  Join each of $p_{i,1}$ and $p_{i,q_i}$ to every
vertex of $T_j$, and join every internal vertex of $P_i$ to every
vertex of $T_{j+1}$.

Now suppose that $i$ is a transition index and $c_i=j$. Join
$p_{i,1}$ to every
vertex of $T_j$; join each of $p_{i,2}$ and $p_{i,q_i}$ to every
vertex of $T_{j+1}$; and join every vertex
\[
 p_{i,3},p_{i,4},\ldots,p_{i,q_i-1}
\]
to every vertex of $T_{j+2}$.  There are no other edges.  In
particular, $T$ is independent, there are no edges between distinct
paths, and
\[
 G-T=P_1+P_2+\cdots+P_{fh}.
\]

For the smallest instance of this construction with $d=3$ and
$k=2$, we have $f=3$, $h=2$, and $n=3[2(2+3)+1]=33.$
In this case $m=0$.  The paths $P_1,P_3,P_5$, for which
$c_i=c_{i+1}$, have order $4$, whereas the transition paths
$P_2,P_4,P_6$ have order $5$.  The resulting graph $G(3,2,33)$ is shown in Figure 3.

\begin{figure}[ht]
\centering
\includegraphics[width=0.75\textwidth]{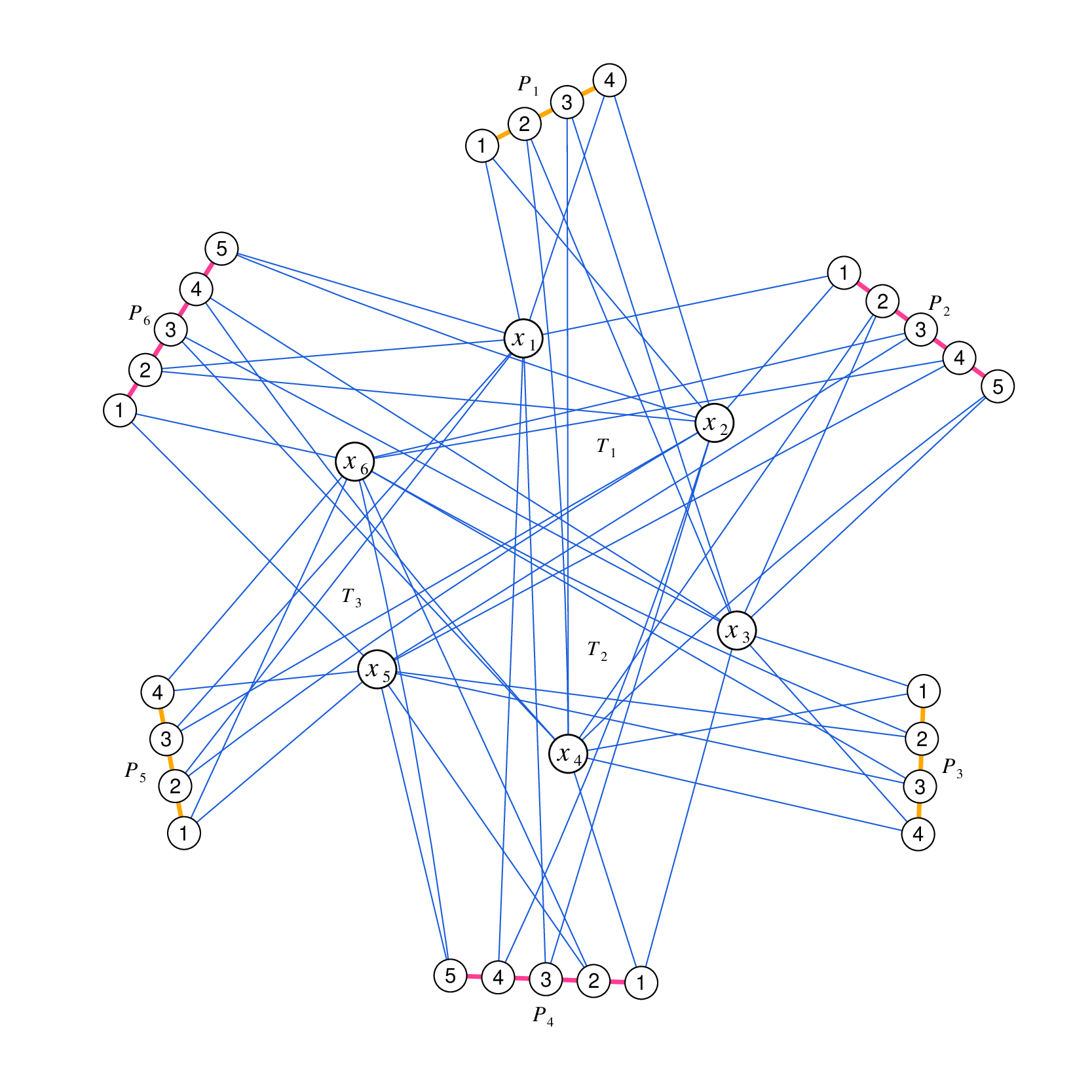}
\caption{The graph $G(3,2,33)$}
\end{figure}

\medskip
\noindent\textit{Verification of the order.}
Before the $m$ additional vertices are added to $P_1$, each of the
$f(h-1)$ paths $P_i$ satisfying $c_i=c_{i+1}$ has order $k+2$, and
each of the $f$ transition paths has order $k+3$. Consequently,
\begin{align*}
 |G|
 &=|T|+f(h-1)(k+2)+f(k+3)+m\\
 &=fh+f(h-1)(k+2)+f(k+3)+m\\
 &=f\bigl(h(k+3)+1\bigr)+m\\
 &=n.
\end{align*}

\medskip
\noindent\textit{Verification of the minimum degree.}
Consider first a path $P_i$ with $c_i=c_{i+1}=j$.  Each endpoint of
$P_i$ has one neighbor on $P_i$
and all $h$ vertices of $T_j$ as neighbors.  Its degree is therefore $1+h=d.$
Each internal vertex has two neighbors on $P_i$ and all $h$ vertices
of $T_{j+1}$ as neighbors, and hence has degree $2+h=d+1.$

For a transition path with $c_i=j$ and with $c_{i+1}$ the cyclic
successor of $c_i$, each of its two endpoints again has one path
neighbor and exactly $h$ neighbors in $T$,
so each endpoint has degree $d$.  Every internal path
vertex has two path neighbors and exactly $h$ neighbors in $T$, so
it has degree $d+1$.

It remains to check the $T$-vertices. Fix $v\in T_j$. For each
vertex $x_i\in T_j$, two distinct path endpoints are adjacent to
$v$: the left endpoint $p_{i,1}$ of $P_i$ and the right endpoint of
the preceding path $P_{i-1}$ in the cyclic sequence; when $i=1$, the
preceding path is $P_{fh}$.  Since $T_j$ contains $h$
vertices, $v$ is adjacent to at least $2h$ path endpoints.  Therefore
\[
 {\rm deg}_G(v)\ge2h=2(d-1)\ge d.
\]
All vertices of $G$ consequently have degree at least $d$, while every
path endpoint has degree exactly $d$.  Hence
\[
 \delta(G)=d.
\]

\medskip
\noindent\textit{Verification of Hamiltonicity.}
For every $i$, the vertex $x_i$ is adjacent to $p_{i,1}$ because
$x_i\in T_{c_i}$, and $x_{i+1}$ is adjacent to $p_{i,q_i}$ because
$x_{i+1}\in T_{c_{i+1}}$.  It
follows that
\[
 x_1P_1x_2P_2\cdots x_{fh}P_{fh}x_1
\]
is a cycle.  It contains every vertex of $T$ and every vertex of every
$P_i$ exactly once, so it is a Hamilton cycle of $G$.

\medskip
\noindent\textit{Exclusion of the forbidden cycle lengths.}
Suppose, for a contradiction, that $G$ has a cycle $C$ of length
$n-s$ for some
\[
 1\le s\le k.
\]
Thus exactly $s$ vertices of $G$ do not belong to $C$.

Every $P_i$ has at least $k+2>s$ vertices.  If $C$ missed an entire
path $P_i$, more than $s$ vertices would be omitted from $C$, which is
impossible.  Hence $C$ meets all $fh$ paths.

Put
\[
 u=|V(C)\cap T|.
\]
We have $u>0$, because otherwise $C$ would be contained in $G-T$, a
disjoint union of paths, which contains no cycle.  Since $T$ is
independent, no two consecutive vertices of $C$ belong to $T$.
Deleting the $u$ vertices of $V(C)\cap T$ from $C$ therefore leaves
exactly $u$ non-null path segments.  Each such segment is contained
in a single component $P_i$ of $G-T$.  Since $C$ meets all $fh$ paths,
there must be at least $fh$ such segments, and hence $u\ge fh.$
On the other hand, $u\le|T|=fh$.  Consequently, $u=fh.$
It follows that $C$ contains every vertex of $T$.  Moreover, deleting
$T$ from $C$ produces exactly $fh$ non-null path segments.  Each
segment lies in one of the $fh$ components $P_1,\ldots,P_{fh}$ of
$G-T$, and $C$ meets every one of these components.  It follows that
exactly one segment lies in each $P_i$.  Because the only edges of
$P_i$ join consecutive vertices in its displayed order, this segment
is a subpath of $P_i$; it consists of consecutive vertices of $P_i$.

For each $i$, denote this \emph{retained subpath} by
\[
 R_i=p_{i,a_i+1}p_{i,a_i+2}\cdots p_{i,q_i-b_i},
 \qquad a_i,b_i\ge0.
\]
More precisely, if $a_i>0$, the \emph{prefix} omitted from $P_i$ is
the initial vertex sequence
\[
 p_{i,1},p_{i,2},\ldots,p_{i,a_i},
\]
which consists of the first $a_i$ vertices of $P_i$; if $a_i=0$, this
prefix is empty.  Similarly, if $b_i>0$, the \emph{suffix} omitted
from $P_i$ is the final vertex sequence
\[
 p_{i,q_i-b_i+1},p_{i,q_i-b_i+2},\ldots,p_{i,q_i},
\]
which consists of the last $b_i$ vertices of $P_i$; if $b_i=0$, this
suffix is empty.  The vertices between these two omitted sequences
form the retained subpath $R_i$.  Since no $T$-vertex is omitted, all $s$
omitted vertices occur in these explicitly described prefixes and
suffixes.  Therefore
\[
 \sum_{i=1}^{fh}(a_i+b_i)=s. \eqno (3)
\]
Since $a_i+b_i\le s$, the retained subpath $R_i$ has at least
\[
 q_i-s\ge(k+2)-k=2
\]
vertices.  In particular, its left and right endpoints are distinct.

We next define the tokens used to record how the cycle $C$ attaches
the retained subpaths to $T$.  Since $R_i$ is a component of $C-T$
and $C$ is a cycle, exactly two edges of $C$ have one endpoint in
$R_i$ and the other in $T$.  Of these two edges, the one incident
with the left endpoint $p_{i,a_i+1}$ of $R_i$ is its \emph{left
token}, and the one incident with the right endpoint
$p_{i,q_i-b_i}$ is its \emph{right token}.  These two endpoints are
distinct, as observed above.  The \emph{actual class} of a token is
the class $T_j$ containing its endpoint in $T$.

We assign a \emph{baseline class} to each token according to the
attachment used by the Hamilton cycle constructed above: the
baseline class of the left token of $R_i$ is $T_{c_i}$, and that of
the right token is $T_{c_{i+1}}$. For each $x_i\in T_j$, the left
token of $R_i$ and the right token of $R_{i-1}$ both have baseline class
$T_j$, where $R_0$ means $R_{fh}$.  Since $|T_j|=h$, exactly $2h$
tokens have baseline class $T_j$.  On the other hand, $C$ contains
every vertex of $T$, and each such vertex has degree $2$ in $C$.
Because $T$ is independent, both incident edges of $C$ are tokens.
Hence exactly
$2h$ tokens have actual class $T_j$ as well.

Orient the cyclic sequence of classes as
\[
 T_1\longrightarrow T_2\longrightarrow\cdots\longrightarrow T_f
 \longrightarrow T_1.
\]
For a token, its \emph{forward movement} is the number of directed
steps in this sequence from its baseline class to its actual class,
using the attachment rules below.  Every forward movement arising in
the construction has length $0$, $1$, or $2$.  In particular, when
$f=2$, a two-step movement
$T_j\longrightarrow T_{j+1}\longrightarrow T_{j+2}=T_j$ is still
recorded as movement $2$, rather than movement $0$.

First consider a path $P_i$ with $c_i=c_{i+1}=j$.  The left token has
baseline class $T_j$.  If $a_i=0$, its endpoint in $R_i$ is
$p_{i,1}$, so its actual class is also $T_j$ and its movement is $0$.
If $a_i\ge1$, that endpoint is $p_{i,a_i+1}$.  It is an internal
vertex of $P_i$: it is not $p_{i,1}$ because $a_i\ge1$, and it is not
$p_{i,q_i}$ because $R_i$ has at least two vertices.  Its actual
class is therefore $T_{j+1}$, so the left token moves one step.  The
right token is analogous: it has movement $0$ if $b_i=0$ and movement
$1$ if $b_i\ge1$.

Next consider a transition path $P_i$ with $c_i=j$. Its left token
has baseline class $T_j$ and behaves
as follows:
\[
 \begin{array}{c|c|c}
 a_i&\text{actual class}&\text{forward movement}\\ \hline
 0&T_j&0\\
 1&T_{j+1}&1\\
 a_i\ge2&T_{j+2}&2.
 \end{array}
\]
Indeed, its endpoint in $R_i$ is respectively $p_{i,1}$,
$p_{i,2}$, or one of $p_{i,3},\ldots,p_{i,q_i-1}$.  In the final case,
the upper bound on the index follows from the fact that at least two
vertices are retained.

The baseline class of the right token of a transition path is
$T_{j+1}$.  If $b_i=0$, its endpoint in $R_i$ is $p_{i,q_i}$, so its
actual class is $T_{j+1}$ and its movement is $0$.  If $b_i\ge1$,
recall that a transition path has $q_i=k+3$ vertices.  Since
$b_i\le s\le k$,
\[
 3\le q_i-b_i\le q_i-1.
\]
The endpoint of the right token in $R_i$ is therefore one of
$p_{i,3},\ldots,p_{i,q_i-1}$, so its actual class is $T_{j+2}$.  The
right token consequently moves one step from $T_{j+1}$ to
$T_{j+2}$.

In every case, the forward movement of a token is at most the number
of vertices omitted at the corresponding end.  Let $D$ be the sum of
the forward movements of all tokens.  Equality (3)
gives
\[
 0\le D\le\sum_{i=1}^{fh}(a_i+b_i)=s\le k. \eqno (4)
\]

In fact, $D>0$.  Indeed, $s>0$ and (3) imply that some $a_i$ or $b_i$ is positive.
For a path satisfying $c_i=c_{i+1}$, any positive $a_i$ or $b_i$
gives movement $1$.  For a transition path, a positive $a_i$ gives
movement $1$ or $2$, and a positive $b_i$ gives movement $1$.
Thus at least one token has positive movement.

We now split each token movement into individual directed steps.  A
traversal of one directed step $T_j\longrightarrow T_{j+1}$ is called
a \emph{unit movement}; here and below class subscripts are read
cyclically, so $T_{f+1}=T_1$.  Thus a token of movement $0$ produces
no unit movement, one of movement $1$ produces one unit movement, and
one of movement $2$ produces two consecutive unit movements.  For
$1\le j\le f$, let $N_j$ be the total number, counted with
multiplicity, of unit movements across the directed step
\[
 T_j\longrightarrow T_{j+1}.
\]
For every class, the number of tokens whose movement starts there is
the number of tokens having that baseline class, namely $2h$; the
number whose movement ends there is the number having that actual
class, also $2h$. Counting token movements with multiplicity therefore
shows that the number of unit movements entering a fixed class equals
the number leaving it. For $T_j$ with $2\le j\le f$, these two numbers are
$N_{j-1}$ and $N_j$; for $T_1$, they are $N_f$ and $N_1$.
Consequently,
\[
 N_f=N_1,
 \qquad
 N_{j-1}=N_j\quad(2\le j\le f). \eqno (5)
\]

Since $D>0$, at least one $N_j$ is positive.  Equalities (5)
then imply that all $N_j$ are equal positive integers.  Since every unit movement contributes $1$ to the
movement of exactly one token, we have
\[
 D=\sum_{j=1}^{f}N_j\ge f=k+1.
\]
This contradicts $D\le k$ from (4).  Hence no cycle of length $n-s$ exists for any $s\in\{1,2,\ldots,k\}$.

We have proved that $G$ has order $n$, is Hamiltonian, has minimum
degree exactly $d$, and contains none of the forbidden cycle lengths. This completes the proof.\hfill $\Box$

\section{Open problems}

Finally, we pose several related open problems.

{\bf Problem 2.} Given an integer $d\ge 3,$ denote by $f(d)$ the least positive integer $p$ such that for every integer $n\ge p,$
there exists a Hamiltonian graph of order $n$ and minimum degree $d$ that contains no $(n-2)$-cycle.

Theorem 1 shows that $f(d)\le 15d-14.$

{\bf Problem 3.} Given an integer $d\ge 3$ and a positive integer $k,$ denote by $g(d,k)$ the least positive integer $q$ such that
for every integer $n\ge q,$ there exists a Hamiltonian graph of order $n$ and minimum degree $d$ that contains no $(n-s)$-cycle for any $s\in\{1,2,\dots,k\}.$

Theorem 2 shows that $g(d,k)\le (k+1)[(d-1)(k+3)+1].$

The following question was posed by the author [9].

{\bf Question 4.} Is it true that every cubic Hamiltonian graph of order $n$ with $n\ge 6$ contains an $(n-2)$-cycle?

We conjecture that the answer to Question 4 is negative. An exhaustive computation shows that every cubic Hamiltonian graph of order $n$ with $6\le n\le 26$ contains an $(n-2)$-cycle.

In 2019, Gir\~{a}o, Kittipassorn and Narayanan posed the following conjecture [3].

{\bf Conjecture 5.} If a graph $G$ of order $n$ and minimum degree at least three contains a Hamilton
cycle, then $G$ contains another cycle of length at least $n-K,$ where $K>0$ is an absolute constant.

Observe that Theorem 2 does not contradict Conjecture 5, since the graph $G(d,k,n)$ constructed in the proof of Theorem 2 has exactly
$[2^{d-2}(d-1)!(d-2)!]^{k+1}$ Hamilton cycles.

\section*{\normalsize Declaration of AI Use}

ChatGPT was used to assist in developing and checking the constructions and proofs. The author independently verified all mathematical arguments, wrote the paper, and takes full responsibility for its content.

\vskip 5mm
{\bf Acknowledgements.} This research was supported by NSFC Grant No. 12271170 and Science and Technology Commission of Shanghai Municipality Grant No. 22DZ2229014.

{\bf Statements and Declarations}

The author has no relevant financial or non-financial interests to disclose.

{\bf Data Availability Statement}

No data were used in this article.
\end{document}